\documentclass[12pt]{article}
\usepackage{amssymb,amsmath,epsfig,amscd,eucal}

\newtheorem{thm}{Theorem}[section]
\newtheorem{defn}[thm]{Definition}
\newtheorem{lem}[thm]{Lemma}
\newtheorem{prop}[thm]{Proposition}

\newtheorem{rem}[thm]{Remark}

\newtheorem{letterthm}{Theorem}
\newtheorem{lettercor}[letterthm]{Corollary}

\newenvironment{pf}{\par\medskip\noindent{\em Proof. }}{\hfill $\square$\par\medskip}
\newenvironment{pfof}[1]{\par\medskip\noindent{\em Proof of #1. }}{\hfill $\square$\par\medskip}

\newcommand{\G}{\mathcal{G}}

\newcommand{\N}{\mathbb{N}}

\newcommand{\Z}{\mathbb{Z}}

\newcommand{\into}{\hookrightarrow}

\input xy
\xyoption{all}

\title{Profinite properties of graph manifolds}
\author{Henry Wilton and Pavel Zalesskii\footnote{Partially supported by CNPq.}}
\date{$18^{\mathrm{th}}$ April 2009}

\begin{document}

\maketitle

\begin{abstract}
Let $M$ be a closed, orientable, irreducible, geometrizable 3-manifold.  We prove that the profinite topology on the fundamental group of $\pi_1(M)$ is efficient with respect to the JSJ decomposition of $M$.  We go on to prove that $\pi_1(M)$ is good, in the sense of Serre, if all the pieces of the JSJ decomposition are.   We also prove that if $M$ is a graph manifold then $\pi_1(M)$ is conjugacy separable.
\end{abstract}

A group $G$ is \emph{conjugacy separable} if every conjugacy class is closed in the profinite topology on $G$.  This can be thought of as a strengthening of residual finiteness (which is equivalent to the trivial subgroup's being closed).  Hempel \cite{hempel_residual_1987} proved that the fundamental group of any geometrizable 3-manifold is residually finite.   In this paper, we investigate which 3-manifolds have conjugacy separable fundamental group.\footnote{A conjugacy separable group has solvable conjugacy problem.  Pr\'{e}aux \cite{pr�aux_conjugacy_2006} has shown that the conjugacy problem is solvable in the fundamental group of an orientable, geometrizable 3-manifold.}  We also study Serre's notion
of goodness, another property related to the profinite topology.

Let $M$ be a compact, connected 3-manifold.    Let $D$ be the closed 3-manifold obtained by doubling $M$ along its boundary. The inclusion $M\into D$ has a natural left inverse.  At the level of fundamental groups it follows that $\pi_1(M)$ injects into $\pi_1(D)$ and that two elements are conjugate in $\pi_1(M)$ if and only if they are conjugate in $\pi_1(D)$.  Hence, if $\pi_1(D)$ is conjugacy separable then so is $\pi_1(M)$.  Therefore, we can assume that $M$ is closed.

Because conjugacy separability is preserved by taking free products \cite{stebe_residual_1970}, we may take $M$ to be irreducible.  As a technical assumption, we shall also assume that $M$ is orientable.\footnote{Conjugacy separability is not a commensurability invariant.  Therefore the orientable case does not immediately imply the non-orientable case.}  Under these hypotheses, $M$ has a canonical \emph{JSJ decomposition}, the pieces of which are either Seifert-fibred or, according to the Geometrization Conjecture, admit finite-volume hyperbolic structures.  By the Seifert--van Kampen Theorem, the JSJ decomposition of $M$ induces a graph-of-groups decomposition of the fundamental group.

Our first theorem asserts that this graph of groups is, from a profinite point of view, well behaved.  If a residually finite group $G$ is the fundamental group of a graph of groups $(\G,\Gamma)$, the profinite topology on $G$ is called \emph{efficient} if the vertex and edge groups of $\G$ are closed and if the profinite topology on $G$ induces the full profinite topologies on the vertex and edge groups of $\G$.

\begin{letterthm}\label{t:Efficiency}
Let $M$ be a closed, orientable, irreducible, geometrizable 3-manifold, and let $(\G,\Gamma)$ be the graph-of-groups decomposition of $\pi_1(M)$ induced by the JSJ decomposition of $M$.  Then the profinite topology on $\pi_1(M)$ is efficient.
\end{letterthm}

Theorem \ref{t:Efficiency} provides the foundation for our main theorems, which relate the profinite completion of $\pi_1(M)$ to the profinite completions of the pieces of the JSJ decomposition.

Our next theorem is a digression regarding goodness, a property introduced by Serre (see I.2.6 Exercise 2 in \cite{serre_galois_1997}).  A group $G$ is \emph{good} if the natural map from $G$ to its profinite completion induces an isomorphism at the level of cohomology with coefficients in any fixed finite $G$-module.  The second author, together with Grunewald and Jaikin-Zapirain, has shown that the fundamental groups of arithmetic hyperbolic 3-manifolds with non-empty boundary are good \cite{grunewald_cohomological_}.

\begin{letterthm}\label{t:Goodness}
Let $M$ be a closed, irreducible, orientable 3-manifold, and suppose that the fundamental groups of the pieces of the JSJ decomposition of $M$ are good.  Then the fundamental group of $M$ is good.
\end{letterthm}

The fundamental group of a Seifert-fibred 3-manifold is good.  We therefore have the following corollary.  Recall that a 3-manifold in which all the pieces of the JSJ decomposition are Seifert-fibred is called a \emph{graph manifold}.

\begin{lettercor}\label{c:Good graph}
If $M$ is a closed, orientable graph manifold then $\pi_1(M)$ is good.
\end{lettercor}

We now return our attention to conjugacy separability.  Martino \cite{martino_proof_2005} has shown that the fundamental groups of Seifert-fibred 3-manifolds are conjugacy separable, and the second author, together with Chagas, has found some examples of hyperbolic 3-manifolds with conjugacy separable fundamental group, including the complement of the figure-eight knot \cite{chagas_figure_2007}.  The fundamental groups of all these 3-manifolds are known to be \emph{subgroup} separable---that is, finitely generated subgroups are closed in the profinite topology \cite{agol_bianchi_2001}.

There are examples of 3-manifolds whose fundamental groups fail to be subgroup separable.  Burns, Karrass and Solitar \cite{burns_note_1987} exhibited a graph manifold whose fundamental group fails to have this property, and Niblo and Wise \cite{niblo_subgroup_2001} showed that the same is true of `most' graph manifolds.  However, our next theorem shows that even these 3-manifolds with poorly-behaved profinite topologies nevertheless are conjugacy separable.

\begin{letterthm}\label{t: Graph manifolds are c.s.}
If $M$ is an orientable graph manifold then $\pi_1(M)$ is conjugacy separable.
\end{letterthm}

This generalizes work of Stebe \cite{stebe_conjugacy_1971}, who showed that the fundamental group of the complement of a hose knot is conjugacy separable.  Note that, for us, graph manifolds include torus bundles over the circle.

Theorem \ref{t: Graph manifolds are c.s.} follows from Theorem \ref{t:Combination theorem}, which provides criteria for graphs of conjugacy separable groups to be conjugacy separable.  The proof of Theorem \ref{t:Combination theorem} proceeds by analysing the action of the profinite completion of $G$ on the profinite Bass--Serre tree.  We are then able to check the hypotheses of Theorem \ref{t:Combination theorem} when the vertex groups are large Seifert-fibred 3-manifolds and to deduce Theorem \ref{t: Graph manifolds are c.s.}.   The hypotheses that Theorem \ref{t:Combination theorem} impose on vertex groups are  more difficult to check for hyperbolic 3-manifolds.  In future work, we intend to check these criteria for hyperbolic 3-manifolds, and to prove that the fundamental group of a closed, orientable, geometrizable 3-manifold is conjugacy separable if the fundamental groups of the pieces of its JSJ decomposition are.

\subsection*{Acknowledgements}

The first author would like to thank Cameron Gordon, Emily Hamilton, Alan Reid and Matthew Stover for useful conversations.

\section{The profinite topology}

\subsection{Basic notions}

Let $G$ be a group.

\begin{defn}
The \emph{profinite topology} on $G$ is the coarsest topology with respect to which every homomorphism from $G$ to a finite group (where the finite group is equipped with the discrete topology) is continuous.
\end{defn}

We are concerned with subsets of $G$ that are closed in the profinite topology.

\begin{defn}
A subset $X$ of $G$ is called \emph{separable} if it is closed in the profinite topology.
\end{defn}

A variety of well known properties can be phrased in terms of separability of certain subsets of $G$.

\begin{defn}
Let $G$ be a group.
\begin{enumerate}
\item If the trivial subgroup $1\subset G$ is separable then $G$ is \emph{residually finite}.
\item If conjugacy classes in $G$ are separable then $G$ is called \emph{conjugacy separable}.
\item If finitely generated subgroups of $G$ are separable then $G$ is called \emph{subgroup separable} or \emph{LERF}.
\end{enumerate}
\end{defn}

\begin{rem}
Note that conjugacy separability, unlike separability of subgroups, is not a commensurability invariant
(\cite{goryaga_example_1986}, \cite{chagas_finite_2007}).\end{rem}

Sometimes a stronger separability property is useful.

\begin{defn}
Let $G$ be a group and $H$ a subgroup.  Then $H$ is \emph{conjugacy distinguished} if, whenever $g\in G$ is not conjugate into $H$, there exists a homomorphism to a finite group $f:G\to Q$ such that $f(g)$ is not conjugate into $f(H)$.
\end{defn}

One approach to the profinite topology on a group is to work with the profinite completion.

\begin{defn}
Consider the inverse system of finite quotients
\[
q:G\to Q
\]
of $G$.  The \emph{profinite completion of $G$} is defined to be the inverse limit
\[
\widehat{G}=\underleftarrow{\lim}\{G\to Q\}.
\]
\end{defn}

We will usually consider residually finite $G$, and in this situation the natural map $G\to\widehat{G}$ is a monomorphism.  As an inverse limit, $\widehat{G}$ comes equipped with a natural topology---the profinite topology on $G$ is the pullback of the topology on $\widehat{G}$.  For a subset $X\subset \widehat{G}$ we write $\overline{X}$ for the closure of $X$ in $\widehat{G}$.

\subsection{Group actions on profinite trees}

For the definitions and basic properties of graphs of groups, and the corresponding group action on a Bass--Serre tree, see \cite{serre_arbres_1977}.  We will be interested in the interaction between graphs of groups and the profinite topology.

\begin{defn}
The profinite topology on the fundamental group of a graph of groups $G=\pi_1(\G,\Gamma)$ is said to be \emph{efficient} if:
\begin{enumerate}
 \item $G$ is residually finite;
 \item it induces the profinite topology on vertex groups $G_v, v\in V(\Gamma)$ and on edge groups $G_e, e\in E(\Gamma)$;
 \item and $G_e, G_v$ are closed in the profinite topology of $G$.
\end{enumerate}
\end{defn}

For example, the fundamental group of a graph of finite groups of bounded order has efficient profinite topology, since it is virtually free (see \cite{dicks_groups_1980}, Chapter IV, Theorem 3.2).  If $\Gamma$ is finite the fact that the profinite topology on $G$ is efficient implies that the profinite completion $\widehat G$ is the profinite fundamental group $\widehat{G}=\pi_1(\widehat\G,\Gamma)$, where $(\widehat\G,\Gamma)$ is the graph of profinite completions of the corresponding vertex and edge groups with the maps $\partial_i$ ($i=0,1$) defined in the obvious (canonical) way. Since a group acts on a tree if and only if it is the fundamental group of a graph of groups (see \cite{dicks_groups_1980}, Theorem 6.1 or \cite{serre_arbres_1977}, Chapter I, Section 5.4), we shall sometimes speak about groups acting on trees with efficient profinite topology. In this case the definition of efficient profinite topology can be given as follows: the profinite topology on $G$ is efficient if $G$ is residually finite, stabilizers of vertices and edges are closed in $G$ and the profinite topology of $G$ induces the (full) profinite topology on the stabilizers.

The idea for establishing the conjugacy separability of the
fundamental group of a graph of groups with efficient profinite
topology is as follows. The profinite completion of such a group
$G=\pi_1(\G,\Gamma)$ is the profinite fundamental group $\widehat
G=\pi_1(\widehat \G,\gamma)$ of the graph $(\widehat \G,\Gamma)$
of profinite completions of the original groups. Therefore we can
consider the standard profinite tree $S(\widehat G)$ associated
with $\widehat G$ and the natural action of $\widehat G$ on it.
This allows us to apply the results of the theory of groups acting
on profinite trees developed in \cite{zalesskii_subgroups_1988},
\cite{zalesskii_geometric_1989},
\cite{zalesskii_fundamental_1989}. We note that $S(\widehat G)$ is
a simply connected profinite graph, which is a stronger property
than being a profinite tree (cf. \cite{zalesskii_geometric_1989}).

The abstract standard tree $S(G)$ associated with $G$ with efficient profinite topology embeds naturally into $S(\widehat G)$ (since all stabilizers are closed) and is dense in it. We distinguish two types of elements of $G$: the elements which stabilize some vertex in $S(G)$ and the elements which move every vertex of $S(G)$. The elements of the second type are called {\it hyperbolic}.

The following result, due to J. Tits (cf. \cite{serre_arbres_1977}), will be used in some of our proofs. We state it here in a form convenient for our purposes.

\begin{prop}\label{4.3.1}
Let $ G=\pi_1(\G,\Gamma)$ be the fundamental group of a graph of groups and assume that $ a\in G $ is hyperbolic. Put
$$ m= \min_{v\in V(S(G))}  \ell[v, av] \ \ \  {\rm and } \ \ \ T_a=  \{ v\in V(S(G)) \mid \ell[v, av]=m\}.$$ Then $ T_a $ is the vertex set of a straight line $($i.e., a doubly infinite chain of $ S(G) )$, that we again denote by $ T_a$, on which  $a$ acts as a translation of amplitude $ m $; furthermore, every $\langle a\rangle $-invariant subtree of  $ S(G)$ contains  $ T_a $, i.e., $T_a$ is unique. Finally if  $ v\in T_a $,  then $ T_a= \langle a\rangle [v, av[ $.

\vskip10pt

\noindent Here, $[v,w]$ denotes the unique geodesic joining vertices $v$ and $w$, and $\ell[v,w]$ its length; also $[v,w[=[v,w] - \{ w\} $. We shall refer to $ T_a $ as the {\it axis} corresponding to the hyperbolic element $ a $.

 $$\cdots \xy (0,0) *{\bullet},(8,2) *{v},(10,0) *{\bullet}
,(20,0) *{\bullet} , (30,0) *{\bullet} ,(38,2) *{av}, (40,0)
*{\bullet} ,(50,0) *{\bullet}
\POS(0,0)\ar(10,0)\POS(10,0)\ar(20,0)\POS(20,0)\ar(30,0) \POS
(30,0)\ar(40,0) \POS(40,0)\ar(50,0)\endxy \cdots$$ \hfill ~$(m=3)$
\end{prop}

The next proposition collects some useful facts about group actions on profinite trees.

\begin{prop}\label{4.3.2}
Let $G=\pi_1(\G,\Gamma)$ be the fundamental group of a finite graph of groups $(\G,\Gamma)$ with efficient profinite topology. Let $b\in G$ be a hyperbolic element and $T_b$ its corresponding axis. Then:
\begin{itemize}
\medskip
\item[\rm (1)] $B=\overline{\langle b\rangle}$ acts freely on the standard graph $S(\widehat\G,\Gamma)$;
\smallskip
\item[\rm (2)]  $B\cong \widehat{\mathbb{Z}}$;
\smallskip
\item[\rm (3)] $\langle b^n\rangle\backslash T_b=\overline{\langle b^n\rangle}\backslash\overline{T_b}$ for every natural number $n$---in particular, the amplitude of the action of $b$ on $T_b$ is the length of the cycle $\langle b^n\rangle\backslash\overline{T_b}$;
\smallskip
\item[\rm (4)] if $\beta\in B$, $v$ is a vertex of $T_b$ and $\beta v\in T_b$, then $\beta\in \langle b\rangle$;
\smallskip
\item[\rm (5)] $T_b$ is a connected component of $\overline {T_b}$ considered as an abstract graph, in other words, the only vertices of $\overline{T_b}$ that are at a finite distance from a vertex of $T_b$ are those of $T_b$;
\smallskip
\item[\rm (6)] $\overline{T_b}$ is the unique minimal profinite $B$-invariant subtree of $S(\widehat G)$ and $\overline{T_b}$ does not contain infinite connected profinite subgraphs.
\end{itemize}
\end{prop}

\begin{pf}
The proofs of (1) and (2) are a word-for-word repetition of the proof of Lemma 4.1 in \cite{ribes_conjugacy_1998}.

To continue the proof of the other items note first that $\overline {T_b} $ is a profinite tree by Theorem 1.15 in \cite{zalesskii_subgroups_1988}.

\medskip\noindent
(3) Since $\overline{\langle b\rangle\ T_b}=B\overline{T_b}$, the tree $T_b$ is a universal covering of $\overline {T_b}/B$ in the abstract (usual) sense. The subgroup of $B$ which leaves $T_b$ invariant is exactly the abstract fundamental group $\pi_1^{abs}(\overline {T_b}/B)$ and $B=\widehat{\pi_1^{abs}}(\overline {T_b}/B)$ (see Lemma 3.5 in \cite{zalesskii_subgroups_1988}). Since $\pi_1^{abs}(\overline {T_b}/B)\cong \mathbb{Z}$, $\langle b\rangle$ is a subgroup of finite index in $\pi_1^{abs}(\overline {T_b}/B)$. But $\langle b\rangle$ is dense in $B$, so it follows that $\langle b\rangle =\pi_1^{abs}(\overline {T_b}/B)$. Hence $\overline {T_b}/B= {T_b}/\langle b\rangle$ and therefore $T_b/\langle b^n\rangle=\overline{T_b}/\overline{\langle b^n\rangle}$ for every natural number $n$.

\medskip\noindent
The proofs of (4) and (5) are the same as in Lemma 4.3 in \cite{ribes_conjugacy_1998}.

\medskip\noindent
(6) Similarly to Example 1.20 in \cite{zalesskii_subgroups_1988} one proves that every connected profinite subgraph of $\overline{T_b}$ is finite. So, if $\Gamma$ is a proper connected $B$-invariant subgraph of $\overline{T_b}$, then $\Gamma$ is finite, which contradicts the freeness of the action of $B$ on $\overline{T_b}$. Then $\overline{T_b}$ is the minimal $B$-invariant subtree, which is unique by Lemma 1.5 in \cite{zalesskii_fundamental_1989}.
\end{pf}

The following lemmas will also be useful.

\begin{lem}\label{4.3.3}
Let $a$ and $b$ be elements of the fundamental group $G=\pi_1(\G,\Gamma)$ of a finite graph of groups $(\G,\Gamma)$ and assume that the profinite topology on $G$ is efficient.  Suppose $a$ and $b$ are conjugate in $\widehat G$ and $a$ is hyperbolic in $G$. Then $b$ is hyperbolic in $G$ and it acts on its corresponding axis $T_b$ with the same amplitude as $a$.
\end{lem}
\begin{pf}
Let $g\in\widehat G$ with $gag^{-1}=b$. Then $b$ acts on $g\overline{T_a}$.  Hence, by Proposition \ref{4.3.2} (6), $g\overline{T_a}=\overline{T_b}$.  Therefore $b$ acts freely on $\overline{T_b}$ and so, by Proposition \ref{4.3.2} (1), $\overline{\langle b\rangle}$ acts freely on $\overline{T_b}$.
Then, by Theorem 2.13 (b) in \cite{zalesskii_subgroups_1988}, $\overline{\langle b\rangle}$ acts freely on
  $S(\widehat G)$ and so on $S(G)$. Thus $b$ is hyperbolic. Moreover, $\overline{\langle a\rangle}\backslash\overline{T_a}\cong\overline{\langle b\rangle}\backslash\overline{T_b}$. Now the result follows from Proposition
\ref{4.3.2} (3).\end{pf}

\begin{lem}\label{4.3.4}
Let $a$ and $b$ be elements of the fundamental group $G=\pi_1(\G,\Gamma)$ of a finite graph of groups $(\G,\Gamma)$ and the  profinite topology on $G$ is efficient.  Suppose $a$ and $b$ are conjugate in $\widehat G$ and $a$ is hyperbolic in $G$. Then for $e\in T_a$ there is a conjugate $b'$ of $b$ in $G$ such that $g_eag_e^{-1}=b'$ for some $g_e\in\widehat G_e$.
\end{lem}

\begin{pf}
Let $g\in\widehat G$ with $gag^{-1}=b$. As was mentioned in the proof of the preceding lemma $g\overline{T_a}=\overline{T_b}$. Then $ge\in \overline{T_b}$.  Choose $b_0\in \overline{\langle b\rangle}$ such that $b_0ge\in T_b$. Since $e$ and $b_0ge$ have the same image in $\widehat G\backslash S(\widehat G)=G\backslash S(G)$ there exists an element $g'\in G$ with $b_0ge=g'e$. Hence $b_0g=g'g_e$ for some $g_e\in \widehat G_e$. Now since $b$ commute with $b_0$ one has $a=g^{-1}b_0^{-1}bb_0g=g_e^{-1}g'^{-1}bg'g_e$. Therefore, putting $b'=g'^{-1}bg'$ we get the result.
\end{pf}

\section{JSJ decompositions of 3-manifolds}

\subsection{Cutting along tori}

The Jaco--Shalen--Johannson decomposition provides a canonical collection of incompressible tori in $M$ \cite{jaco_seifert_1979,johannson_homotopy_1979}.   For an accessible account of this theorem, see \cite{neumann_canonical_1997}.

\begin{thm}
Any closed, orientable, irreducible 3-manifold $M$ contains a finite embedded collection of disjoint, incompressible tori $\mathcal{T}=\bigcup_i T_i$ such that if  $\{M_j\}$ are the connected components of $M\smallsetminus \mathcal{T}$, each $M_j$ is either a Seifert-fibred manifold or is atoroidal.  A minimal such collection is unique up to isotopy.
\end{thm}

The minimal such decomposition of $M$ is called the \emph{JSJ decomposition of $M$}, and induces a graph-of-groups decomposition $(\G,\Gamma)$ of the fundamental group $G$.  We will call this the \emph{JSJ decomposition of $G$}.

Because the tori are incompressible, each piece $M_j$ of the resulting decomposition has incompressible, toral boundary.  We shall call each component of the boundary a \emph{cusp} of $M_j$.  A subgroup $H$ of
$\pi_1(M_j)$ is \emph{peripheral} if it is conjugate to the fundamental group of a cusp.

The theorem asserts that each piece is either Seifert-fibred or atoroidal.  See subsection \ref{ss:Seifert-fibred} for more on Seifert-fibred manifolds.  A 3-manifold $M$ is called \emph{atoroidal} if every $\pi_1$-injective map of a torus into $M$ is homotopic into the boundary.  Thurston's famous Geometrization Conjecture implies that each atoroidal piece admits a hyperbolic structure of finite volume \cite{thurston_three-dimensional_1982}.  We call such a 3-manifold \emph{geometrizable}.  Perelman has announced a proof of the Geometrization Conjecture \cite{perelman_entropy_2002,perelman_ricci_2003,perelman_finite_2003}.  In the case when every piece of the JSJ decomposition is Seifert-fibred, $M$ is called a \emph{graph manifold}.

As all our theorems are immediate when the JSJ decomposition of $M$ is trivial, we shall usually assume that the JSJ decomposition is non-trivial---that is, the induced graph of groups $(\G,\Gamma)$ has at least one edge.

\emph{Torus bundles} are a special case of graph manifolds.  These can be constructed from a torus crossed with an interval by identifying the two boundary components using an automorphism of the torus.  If the automorphism is hyperbolic---that is, it corresponds to an element of $GL_2(\Z)$ with distinct real eigenvalues---then the resulting manifold admits a geometric structure modelled on the Lie group $\mathrm{Sol}$.  Otherwise, the resulting bundle is Seifert-fibred.

\subsection{Seifert-fibred manifolds}\label{ss:Seifert-fibred}

We will not give a precise definition of Seifert-fibred 3-manifolds here.  They can be defined as `Seifert bundles' over cone-type 2-orbifolds.  The fibres of the bundle are circles.   For our purposes, we just need to know that the fundamental group $G$ of a compact Seifert-fibred 3-manifold $M$ fits into a short exact sequence
\[
1\to Z\to G\stackrel{p}{\to}\pi_1(O)\to 1
\]
where $Z=\langle z\rangle$ is cyclic and $O$ is a cone-type 2-orbifold.  For more on Seifert-fibred manifolds see, for instance, \cite{scott_geometries_1983}.

Seifert-fibred manifolds have well-behaved profinite topologies---they are double-coset separable \cite{niblo_separability_1992} and, more importantly for our purposes, conjugacy separable.

\begin{thm}[Martino, \cite{martino_proof_2005}]\label{t: Martino}
If $M$ is a Seifert-fibred manifold then $\pi_1(M)$ is conjugacy separable.
\end{thm}

We will be interested in gluing Seifert-fibred spaces along cusps, corresponding to (conjugacy classes of) subgroups of the form $\langle\delta\rangle\oplus Z\cong\Z^2$ where $\delta$ maps to a boundary component of $O$.  In particular, for us $Z$ will always be infinite.

Any Seifert-fibred manifold with a cusp is particularly simple---its Seifert bundle structure is \emph{virtually trivial}.  Let $M$ be a Seifert-fibred manifold with a cusp, and $O$ the base orbifold.  By standard theory, $O$ has a finite-sheeted cover $\Sigma$ that is a genuine surface.  Pulling back gives a finite-sheeted cover $M'\to M$ fitting into the short exact sequence
\[
1\to Z\to\pi_1(M')\to\pi_1(\Sigma)\to 1.
\]
Since $M$ has a cusp $\Sigma$ has a boundary component, so $\pi_1(\Sigma)$ is free and the extension splits:
\[
\pi_1(M')\cong Z\rtimes\pi_1(\Sigma).
\]
If $\Sigma$ is a surface of negative Euler characteristic then the
orbifold $O$ is called \emph{hyperbolic} and $M$ is called
\emph{large}.  (More generally, a 3-manifold is large if its
fundamental group has a finite-index subgroup that surjects onto a
non-abelian free group.)

If a Seifert-fibred manifold $M$ is not large but does have an incompressible toral boundary component then the base orbifold $O$ is Euclidean, and is either an annulus, a M\"obius band or a disc with two cone points of order two.  If $M$ is orientable and the base orbifold is an annulus then $M$ is homeomorphic to a direct product of a torus and an interval; as the two boundary components are parallel, such a Seifert-fibred piece can only occur in the JSJ decomposition of a closed 3-manifold if the two boundary components are identified, so if the closed manifold is a torus bundle over a circle.

The remaining two Euclidean orbifolds, the M\"obius band and the disc with two cone points of order two, arise as the base orbifolds of two different Seifert-fibred structures on the twisted interval bundle over the Klein bottle.  If $M$ is homeomorphic to this interval bundle then $\pi_1(M)$ is isomorphic to $\Z\rtimes\Z$ and the index-two abelian subgroup corresponds to a two-sheeted covering space homeomorphic to the trivial interval bundle over the torus.

\subsection{Acylindrical splittings}

Let $k$ be a positive integer.  A graph of groups $\G$ is called \emph{$k$-acylindrical} if, for any $g\in G\smallsetminus 1$, the diameter of the fixed point set of $g$ in the Bass--Serre tree is at most $k$.

\begin{lem}\label{l: 2-acylindrical lemma}
Let $M$ be a closed, orientable, irreducible 3-manifold. If every vertex space of the JSJ decomposition $\G$ is either hyperbolic or a large Seifert-fibred space then $\G$ is $2$-acylindrical.
\end{lem}
\begin{pf}
Let $G=\pi_1(M)$ and let $T$ be the Bass--Serre tree corresponding to $\G$.  Consider an arc $I$ of length 3 in $T$, consisting of consecutive edges $\hat{e}_0$, $\hat{e}_1$ and $\hat{e}_2$, separated by vertices $\hat{v}_1$ and $\hat{v}_2$.  For $i=1,2,3$ let $e_i$ be the edge of $\G$ covered by $\hat{e}_i$ and for $i=1,2$ let $v_i$ be the vertex of $\G$ covered by $\hat{v}_i$.

Suppose that $v_1$ corresponds to a hyperbolic piece of $M$.  Then, for any $g\in G_{v_1}$, $G_{e_0}^g\cap G_{e_1}=1$ and it follows that the intersection of the stabilizers of $\hat{e}_0$ and $\hat{e}_1$ are trivial, and therefore the stabilizer of $I$ is trivial.  The same holds if $v_2$ corresponds to a hyperbolic piece.  We can therefore assume that $v_1$ and $v_2$ correspond to Seifert-fibred pieces $M_1$ and $M_2$ respectively.  For each $i$, let $G_i=\pi_1(M_i)$, let $O_i$ be the underlying orbifold of $M_i$ and let $Z_i$ be the cyclic subgroup of $G_i$ that corresponds to a regular fibre.  Without loss of generality, we can assume that $G_i$ stabilizes $\hat{v}_i$.

By hypothesis, each $O_i$ is hyperbolic.  It follows that any two peripheral subgroups of $O_i$ have trivial intersection, and hence that $G_{\hat{e}_0}\cap G_{\hat{e}_1}=Z_1$ and, likewise, $G_{\hat{e}_1}\cap G_{\hat{e}_2}=Z_2$.  If $I$ has a non-trivial stabilizer then it follows that $Z_1$ and $Z_2$ have non-trivial intersection, and hence are equal, as they are both direct factors in $G_{e_1}$.  But if this is the case then the Seifert-fibred structures of $M_1$ and $M_2$ coincide on their common torus.  This contradicts the minimality of the JSJ decomposition of $M$.
\end{pf}

\begin{lem}\label{l: 4-acylindrical lemma}
Let $M$ be a closed, orientable, irreducible 3-manifold. Either $M$ has a finite-sheeted covering space that is a torus bundle over a circle or the JSJ decomposition $\G$ is $4$-acylindrical.
\end{lem}
\begin{pf}
By Lemma \ref{l: 2-acylindrical lemma}, it remains to deal with the case when some Seifert-fibred vertices $\{M_i\}$ of $\G$ are not large.  Let $\{N_j\}$ be the remaining, hyperbolic and large Seifert-fibred, pieces.  If $M$ is not a torus bundle over a circle then every $M_i$ is homeomorphic to a twisted interval bundle over a Klein bottle.  Unless $M$ is finitely covered by a torus bundle, no two $M_i$ are adjacent.  There is a two-sheeted covering space $M'$ constructed by taking two copies of every $N_j$, replacing every $M_i$ with its double cover $M'_i$ that is homeomorphic to the trivial interval bundle over the torus, and gluing appropriately.  Equivalently, $M'$ can be seen as corresponding to the kernel of the map $\pi_1(M)\to\Z/2$ that assigns to each element of $\pi_1(M)$ the parity of the sum of its intersection number with the core Klein bottles of the $M_i$.  No two $M'_i$ are adjacent in $M'$, and the JSJ decomposition of $M'$ is obtained by collapsing the $M'_i$ to tori.

Let $T$ be the Bass--Serre tree of $T$ and let $\G'$ be the decomposition of $G'=\pi_1(M')$ induced by the action of $G'$ on $T$.  By the above discussion, $\G'$ is obtained from the JSJ decomposition of $M'$ by subdividing certain edges once.  Therefore, if $g\in G'$ then the fixed point set of $g$ in $T$ has diameter at most 4, by Lemma \ref{l: 2-acylindrical lemma}.

It remains to consider $g\in G\smallsetminus G'$.  Suppose such a $g$ acts elliptically on $T$.  Then $g$ stabilizes a vertex corresponding to a twisted interval bundle over the circle, with stabilizer isomorphic to $\Z\rtimes\Z$.  The stabilizers of the two incident edges are both equal to the abelian subgroup of index two.  This subgroup is contained in $G'$, and $g$ interchanges the two incident edges, so the fixed point set of $g$ has diameter equal to $0$.
\end{pf}

\section{Efficiency}

In this section, we make use of Hamilton's results in \cite{hamilton_abelian_2001} to show that the profinite topology on the fundamental group of a geometrizable 3-manifold equipped with its JSJ decomposition is efficient.

\subsection{Normal $n$-characteristic subgroups}

Following Hamilton \cite{hamilton_abelian_2001}, we define $n$-characteristic subgroups of the fundamental group of a 3-manifold with incompressible toral boundary.  Hamilton's definition is phrased in terms of covering spaces.  We state it here in terms of the fundamental group.

\begin{defn}[Hamilton, \cite{hamilton_abelian_2001}]
Let $M$ be a compact 3-manifold with incompressible toral boundary.  Let $G=\pi_1(M)$ and let $n$ be a positive integer.  A finite-index subgroup $K$ is called $n$-characteristic if
\[
A\cap K=nA
\]
for any peripheral subgroup $A$ of $G$.
\end{defn}

Note that $n$-characteristic subgroups need not be characteristic in the usual sense of invariant under group automorphisms, or even normal.  To prove that the profinite topology on $\pi_1(M)$ is efficient, we will need to construct a rich supply of normal, $m$-characteristic subgroups of the pieces of the JSJ decomposition.  Indeed, we will use the following result.

\begin{thm}\label{t:Normal characterisic subgroups}
Let $M$ be a closed, orientable, irreducible, geometrizable 3-manifold.  Let $\{G_i\}$ be the fundamental groups of the pieces of the JSJ decomposition of $M$. For every integer $n\in\N$ there exists an integer $\nu_n$ such that for every group $G_i$ there exists a normal $\nu_n n$-characteristic subgroup $K_{i,n}$ of $G_i$.
\end{thm}

Hamilton finds normal $n$-characteristic subgroups for the hyperbolic pieces.

\begin{lem}[Hamilton, Lemma 5 of \cite{hamilton_abelian_2001}]\label{l: Characteristic subgroups of hyperbolic pieces}
Let $M_1,\ldots,M_n$ be a finite collection of finite-volume hyperbolic 3-manifolds and let $m$ be a positive integer.  Then there exists an integer $\nu$ and normal $\nu m$-characteristic subgroups $K_{i,m}\lhd\pi_1(M_i)$ for each $i$.
\end{lem}

Hamilton also finds $n$-characteristic subgroups of Seifert-fibred 3-manifolds, although they are not normal.

\begin{lem}[Hamilton, Lemma 6 of \cite{hamilton_abelian_2001}]\label{l: Hamilton sf}
Let $M$ be a compact, orientable Seifert-fibred manifold with non-empty incompressible boundary and let $G=\pi_1(M)$.  There exists a positive integer $\mu_M$ such that, for any $n\in\N$, there exists a $\mu_M n$-characteristic subgroup $L_n\subset G$.
\end{lem}

This is easily adapted to produce normal $m$-characteristic subgroups.

\begin{lem}\label{l: Characteristic subgroups of Seifert-fibred pieces}
Let $M$ be a compact, orientable Seifert-fibred manifold with non-empty incompressible boundary and let $G=\pi_1(M)$.  There exists a positive integer $\mu_M$ such that, for any $n\in\N$, there exists a normal $\mu_M n$-characteristic subgroup $K_n\lhd G$.
\end{lem}
\begin{pf}
Let $n$ be a  positive integer, and let $L_n$ be the $\mu_M n$-characteristic subgroup provided by Lemma \ref{l: Hamilton sf}.  Set
\[
K_n=\bigcap_{g\in G}gL_ng^{-1}.
\]
This is a finite-index, normal subgroup of $G$.  It remains to see that $K_n$ is indeed $\mu_M n$-characteristic.

Let $A$ be a peripheral subgroup.  For any $g\in G$, $A'=g^{-1}Ag$ is also a peripheral subgroup, so by the construction of $L_n$ we have
\[
g^{-1}Ag\cap L_n=\mu_M ng^{-1}Ag.
\]
Conjugating by $g$, we see that
\[
A\cap gL_ng^{-1}=\mu_M nA
\]
and so $A\cap K_n=\mu_M nA$ as required.
\end{pf}

Theorem \ref{t:Normal characterisic subgroups} follows easily from Lemmas \ref{l: Characteristic subgroups of hyperbolic pieces} and \ref{l: Characteristic subgroups of Seifert-fibred pieces}.

\begin{pfof}{Theorem \ref{t:Normal characterisic subgroups}}
Let $G_i$ be the fundamental groups of the hyperbolic pieces of the JSJ decomposition and let $H_j$ be fundamental groups of the Seifert-fibred pieces.  For each $j$, Lemma \ref{l: Characteristic subgroups of Seifert-fibred pieces} provides a positive integer $\mu_j$. Let $m=n\prod_j \mu_j$.  By Lemma \ref{l: Characteristic subgroups of hyperbolic pieces}, there exists $\nu_m$ such that every $G_i$ admits a normal $\nu_m m$-characteristic subgroup.  Now every $H_j$ admits a normal $\nu_m m$-characteristic subgroup, as well.
\end{pfof}

\subsection{Proof of Theorem \ref{t:Efficiency}}

We will need the following theorem of Long and Niblo.

\begin{thm}[Long, Niblo \cite{long_subgroup_1991}]\label{t: Peripheral subgroups are separable}
Let $M$ be a compact 3-manifold and let $S$ be an incompressible component of the boundary of $M$.  Let $H$ be a subgroup of $G=\pi_1(M)$ conjugate to $\pi_1(S)$.  Then $H$ is separable in $G$.
\end{thm}

We are now in a position to prove Theorem \ref{t:Efficiency}.

\begin{thm}\label{efficient}
Let $M$ be a closed, orientable, irreducible, geometrizable 3-manifold.  The profinite topology on $G=\pi_1(M)$ equipped with the graph of groups induced by the JSJ decomposition of $M$ is efficient.
\end{thm}
\begin{pf}
First, by \cite{hempel_residual_1987}, $G$ is residually finite.

Observe that Theorem \ref{t:Normal characterisic subgroups} implies that every finite-index subgroup of an edge group contains a finite-index subgroup that extends to a finite-index subgroup of any adjacent vertex group.  So vertex groups induce the full profinite topology on edge groups.

Next we shall see that $G$ induces the full profinite topology on the vertex groups.  Let $v$ be a vertex and let $N_v$ be a finite-index normal subgroup of $G_v$. It suffices to find a finite-index subgroup $N$ of $G$ such that $N\cap G_v\subset N_v$.  There exists an integer $n$ such that, for any edge $e$ incident at $v$, $nG_e\subset N_v\cap G_e$.  By Theorem \ref{t:Normal characterisic subgroups} there exists $\nu_n$ such that for every vertex $u$ there exists a finite-index normal subgroup $M_{u,n}\lhd G_u$ such that for any edge $f$ incident at $u$, $M_{u,n}\cap G_f=\nu_n nG_f$.  In the case $u=v$, by replacing $M_{v,n}$ with $M_{v,n}\cap N_v$ we may assume that $N_v$ contains $M_{v,n}$.  Let $\phi:G\to G_n$ be the quotient obtained by quotienting each vertex $G_v$ by $M_{v,n}$.  Then $G_n$ is a virtually free group, hence residually finite, and $\phi(G_v)$ is finite, so there exists a finite-index, normal subgroup $K\lhd G_n$ such that $\phi(G_v)\cap K=1$.  Let $N=\phi^{-1}(K)$.  Then
\[
N\cap G_v=M_{v,n}\subset N_v
\]
as required.

As edge groups are separable in vertex groups by Theorem \ref{t: Peripheral subgroups are separable},
it remains to prove that any vertex group $G_v$ is separable in $G$.  Let $g\in G\smallsetminus G_v$.
Then the reduced form
\[
g=g_0e_1g_1\ldots e_kg_k
\]
in $\pi_1(\G,v)$ is non-trivial (see 5.2 in \cite{serre_arbres_1977}).  For each vertex $u$, let $K_u$ be a finite-index normal subgroup of $G_u$ such that, whenever $g_i\in G_{u}$, $g_i\notin G_{e_{i-1}}K_u$ and $g_i\notin G_{e_i}K_u$.  (The existence of such a $K_u$ is guaranteed by Theorem \ref{t: Peripheral subgroups are separable}.)  There is an $n$ such that, for each vertex $u$ and each edge $e$ incident at $u$, $nG_e\subset K_u\cap G_e$.  By Theorem \ref{t:Normal characterisic subgroups}, there is $\nu_n$ such that, for each vertex $u$ there is a finite-index normal subgroup $L_{u,n}\lhd G_u$ such that $L_{u,n}\cap G_e=\nu_nnG_e$ for each incident edge $e$.  Let $M_{u,n}=K_u\cap L_{u,n}$, and as before let $G_n$ be the fundamental group of the graph of finite groups obtained by quotienting each $G_u$ by $M_{u,n}$, and let $\phi:G\to G_n$ be the natural map.  Then \[
\phi(g)=\phi(g_0)e_1\phi(g_1)\ldots e_k\phi(g_k)
\]
is a (non-trivial) reduced form for $\phi(g)$, so $\phi(g)$ is not in $\phi(G_u)$.  Once again, $G_n$ is virtually free and so residually finite, and it follows that there is a finite quotient $\psi:G\to Q$ such that $\psi(g)\notin \psi(G_v)$, as required.
\end{pf}

\section{Goodness}

We begin with the definition introduced by Serre (see I.2.6 Exercise 2 in \cite{serre_galois_1997}).

\begin{defn}
Let $G$ be a group and $\widehat G$ its profinite completion.  The group $G$ is called \emph{good} if the
homomorphism at the level of cohomology
\[
H^n(\widehat{G},A)\to H^n(G,A)
\]
is an isomorphism for every finite $G$-module $A$.
\end{defn}

Here the cohomology of $\widehat G$ is the continuous cohomology and cohomolgy of $G$ is usual one.

\begin{prop}\label{p: Good Seifert}
If $M$ is a compact Seifert-fibred manifold then $G=\pi_1(M)$ is good.
\end{prop}

\begin{pf} Since $G/Z$ is good the result follows from   I.2.6 Exercise 2 (c) in \cite{serre_galois_1997}\end{pf}

Theorem \ref{t:Goodness} follows immediately from Theorem \ref{t:Efficiency} and the next proposition. Applying Proposition \ref{p: Good Seifert}, Corollary \ref{c:Good graph} is immediate.

\begin{prop}\label{techprop}
Let $G$ be the fundamental group of a finite graph of good groups and suppose the profinite topology on $G$ is efficient. Then $G$ is good.\end{prop}
\begin{pf}
If $\Gamma$ is finite the fact that the profinite topology on $G$ is efficient implies that the profinite
completion $\widehat G$ is the profinite fundamental group $G=\pi_1(\widehat\G,\Gamma)$,
where $(\widehat\G,\Gamma)$ is the graph of profinite completions of the corresponding vertex and edge groups
with the maps $\partial_i$ ($i=0,1$) defined in the obvious (canonical) way. Therefore we can consider the
standard profinite tree $S(\widehat G)$ associated with $\widehat G$ and the natural action of $\widehat G$ on
it. We note that the abstract standard tree $S(G)$ embeds naturally into $S(\widehat G)$ (since all stabilizes
are closed) and is dense in it. Thus
we have the following commutative diagram (cf. (1.15), (3.7), (3.8) in \cite{zalesskii_subgroups_1988}) whose rows are short exact sequences associated to trees $S(G)$ and $S(\widehat G)$:

$$\xymatrix{0\ar[r]&\bigoplus_{e\in E(\Gamma)}\widehat\Z[[\widehat G
/\widehat G_e]]\ar[r]&\bigoplus_{v\in
V(\Gamma)}\widehat\Z[[\widehat G/\widehat G_v]]\ar[r]&\widehat\Z\ar[r]&0\\
0\ar[r]&\bigoplus_{e\in
E(\Gamma)}\Z[G/G_e]\ar[r]\ar[u]&\bigoplus_{v\in
V(\Gamma)}\Z[G/G_v]\ar[r]\ar[u]&\Z\ar[r]\ar[u]&0}$$

Applying $Hom_{\Z[G]}(-,M)$ to the second row and continuous $Hom_{\widehat \Z[[\widehat G]]}(-,M)$ to
the first row, where $M$ is a finite  $G$-module we get the commutative diagram of the Mayer-Vietoris sequence associated to $G$ and $\widehat G$:
$$\begin{matrix}
\prod_{v\in V(\Gamma)}H^{n-1}(G_v,M)& \rightarrow & H^n(G,M)& \rightarrow & \prod_{e\in E(\Gamma)}H^n(G_e,M)&\rightarrow&\cdots\\
 \uparrow& & \uparrow & &\uparrow \\
\prod_{v\in V(\Gamma)}H^{n-1}(\widehat{G}_v,M)& \rightarrow & H^n(\widehat G,M)& \rightarrow&\prod_{e\in E(\Gamma)}H^n(\widehat G_v,M)& \rightarrow &\cdots
\end{matrix}$$
where the vertical maps are induced by the natural embedding of the groups into their profinite completions. Since $G_e$ and $G_v$ are good the left vertical map and the right vertical map are isomorphisms, so the middle vertical map is an isomorphism as well. Since $H^0(G,M)=M^G=M^{\widehat G}=H^0(\widehat G,M)$ the result follows.
\end{pf}

\section{Profinite completions of groups acting on trees}

\subsection{A combination theorem for conjugacy separable groups}


We shall apply the technology of profinite group actions on profinite trees to prove a combination theorem for conjugacy separable groups.  Theorem \ref{t: Graph manifolds are c.s.} is an application of this.

\begin{defn}
Let $(\G,\Gamma)$ be graph of groups and suppose that the
profinite topology on $G=\pi_1(\G,\Gamma)$ is efficient.  Then
$(\G,\Gamma)$ is \emph{profinitely $k$-acylindrical} if the
corresponding action of the profinite completion $\widehat{G}$ on
the profinite Bass--Serre tree $S(\widehat{G})$ is
$k$-acylindrical---that is, whenever
$\gamma\in\widehat{G}\smallsetminus 1$, the fixed point set of
$\gamma$ in $S(\widehat{G})$ has diameter at most $k$. Note that
by Corollary 4 in \cite{HZ} this means that any element $1\neq
g\in \widehat G$ can fix at most $k$ edges  in any (profinite)
geodesic $[v,w]$ of $S(\widehat{G})$.
\end{defn}

We can now state our main technical result.

\begin{thm}\label{t:Combination theorem}
Let $(\G,\Gamma)$ be a finite graph of groups with conjugacy separable vertex groups.  Let $G=\pi_1(\G,\Gamma)$, and suppose that the profinite topology on $G$ is efficient and that $(\G,\Gamma)$ is profinitely 2-acylindrical.  For any vertex $v$ of $\Gamma$ and incident edges $e$ and $f$, suppose furthermore that the following conditions hold:
\begin{enumerate}
\item for any $g\in G_v$ the double coset $G_egG_f$ is separable in $G_v$;
\item the edge group $G_e$ is conjugacy distinguished in $G_v$;
\item the intersection of the closures of $G_e$ and $G_f$ in the profinite completion of $G_v$ is equal to
the profinite completion of their intersection, i.e. $\overline{G}_e\cap\overline{G}_f=\widehat{G_e\cap G_f}$.
\end{enumerate}
Then $G$ is conjugacy separable.
\end{thm}
\begin{pf}
As the profinite topology on $G$ is efficient, $S=S(G)$ is embedded in $S(\widehat G)$. Fix a connected
transversal $\Sigma$ in $S$. Let $a,b\in G$, and assume that $\gamma a\gamma^{-1}=b$ for some $\gamma\in\widehat G$. Our aim is to show that $gag^{-1}=b$ for some $g\in G$.

\medskip
{\it Case 1.} One of the elements $a,b$ is conjugate to an element of a vertex group $G_v$ (i.e. $a$ is not hyperbolic). Then by Lemma \ref{4.3.3}   $b$ is not hyperbolic as well. Thus we may assume that $a\in G_v, b\in G_w$ for some vertices $v,w\in \Sigma$.  If $\gamma$ belongs to $\widehat G_v$, then the result follows from the efficiency of the profinite topology and conjugacy separability of $G_v$. Otherwise, by Theorem 3.12 of \cite{zalesskii_subgroups_1988} we have $a\in\alpha\overline G_e\alpha^{-1}$, $b\in\beta\overline G_{e'}\beta^{-1}$ for some $\alpha\in \widehat G_v$, $\beta\in\widehat G_w$ and some edges $e,e'$ of $S(\widehat G)$. It follows that $a$ and $b$ are conjugate into the profinite completions of some edge groups in $\widehat G$.  Then $a$ and $b$ are conjugate in $G$ to elements of edge groups because edge groups are conjugacy
distinguished.  So we may assume that $a\in G_e$, $b\in G_{\bar e}$ for some $e,\bar e\in \Sigma$. Then $b$ fixes $[\bar e,\gamma  e]$. As $(\G,\Gamma)$ is profinitely 2-acylindrical, this geodesic can have at most two edges $\bar e$ and $\gamma e$. Let $u$ be their common vertex. Since $u\in \Sigma$, $\gamma$ has to be in $\widehat G_u$ the case already considered above. This completes the proof of Case 1.

\medskip
{\it Case 2}. The element $a$ is hyperbolic. By Lemma \ref{4.3.3}, $b$ is hyperbolic and  acts freely on $ S(G)$. Consider the axis  $ T_a  $ and  $ T_b  $  corresponding to $ a $ and $b$ (see Proposition \ref{4.3.1}). By Lemma \ref{4.3.3} they act with the same amplitude $m$.

By Lemma \ref{4.3.4} we can assume  that  $\gamma$ is in $\widehat G_e$ for some $e\in T_a$. Then $e\in T_b$, by Proposition \ref{4.3.1}. We need to arrange that $\gamma$ fixes longer geodesics in $T_a$ than the geodesic whose only edge is $e$. Suppose that $P$ is a finite geodesic in $T_a$ which has $e$ as one of its edges and such that $\gamma\in \widehat G_P$ (where $\widehat G_P$ means the intersection of the stabilizers in $\widehat G$ of all edges of $P$);  we shall show that $\gamma$ can be replaced by an element that lies in the closure of the intersection of the edge stabilizers of a geodesic strictly containing $P$. Note first that because $(\G,\Gamma)$ is profinitely 2-acylindrical, $P$ has  at most two edges $e$ and $\bar e$ and by the third hypothesis in this case $\gamma\in \widehat G_e\cap \widehat G_{\bar e}=\widehat{G_e\cap G_{\bar e}}$.  Thus $\widehat G_P=\overline{G}_P$.

Let $e_1$ be an edge of $T_a\smallsetminus P$ connected to $P$, write $v$ for the common vertex of $e_1$ and
$P$, and write $P_+$ for the graph with edges those of $P$ together with $e_1$. Let $e_2=\gamma e_1\in\overline T_b$.

First note that $e_2\in T_b$. Indeed, if $e'$ is an edge of $T_b$ there is a path in $S(G)$ connecting $e'$ to $e_1$, and so since $e_1,e_2$ share a vertex there is path connecting $e'$ to $e_2$.  However, since $e_2=\gamma e_1\in\gamma\overline{T_a}=\overline{T_b}$ and $e'\in\overline{T_b}$, it follows that the geodesic $[e',e_2]$ lies in $\overline{T_b}$. The abstract connected component of $\overline{T_b}$ containing $e'$ is precisely $T_b$ (by Proposition
 \ref{4.3.2} (5)), and so we conclude that $e_2\in T_b$.

Now because $S(G)/G=S(\widehat G)/\widehat G$ we have $ge_1=e_2$ for some $g\in G$ and since $v$ is a common vertex of $e_1$ and $e_2$ the element $g\in G_v$. Since $e_1=g^{-1}e_2=\gamma^{-1}e_2$ the element $\gamma_1=\gamma g^{-1}$ is in $\overline{G_{e_2}}$.  Moreover, $\overline{G_{e_2}G_P}\cap G_v=G_{e_2}G_P$ as double cosets of edge groups are separable. Therefore, because $g=\gamma_1^{-1}\gamma$ we can find $h_1\in G_P$, $h_2\in G_{e_2}$ with $g=h_2h_1$. Set $\gamma_+=h_1^{-1}\gamma$. Thus

$$\gamma_+e_1=h_1^{-1}e_2=g^{-1}h_2e_2=g^{-1}e_2=e_1,$$ and so $\gamma_+\in\overline{G_{e_1}}$.
We also have $\gamma_+\in \overline{G_P}$. Then,  by assumption 3 if $P$ has just one edge and by profinite 2-acylindricity if $P$ has two edges,
$\gamma_+\in \widehat G_{P_+}$. We may therefore replace $\gamma$ by $\gamma_+$ and $b$ by $h_1^{-1}bh_1$ and so assume that $\gamma$ is  in $\widehat G_{P_+}$.

But, because $(\G,\Gamma)$ is profinitely 2-acylindrical, a path having more than three edges can only have trivial stabilizers. This finishes the
proof.
\end{pf}

\subsection{Graph manifolds}

We shall see that the hypotheses of Theorem \ref{t:Combination theorem} apply to the JSJ decomposition of a graph manifold when every Seifert-fibred piece is large.  Therefore, the fundamental group of such a graph manifold is conjugacy separable, and we shall deduce the general case.

Theorem \ref{t:Efficiency} asserts that the JSJ decomposition is efficient.  As Seifert-fibred 3-manifolds have double-coset separable fundamental group \cite{niblo_separability_1992}, it remains only to check hypotheses 2 and 3, and that the JSJ decomposition is profinitely 2-acylindrical.

\begin{lem}\label{conjugacy distinguished}
Let $G$ be the fundamental group of a Seifert-fibred space. Then every peripheral subgroup is conjugacy distinguished.
\end{lem}
\begin{pf}
Suppose $g$ is not conjugate to any element of $A$. Since $Z\leq A$ the image of $g$ in $G/Z$ is not conjugate to an element of $A/Z$. But $G/Z$ is virtually free, so by \cite{ribes_conjugacy_1998} $A/Z$ is conjugacy distinguished. Hence there is a finite quotient of $G/Z$ where the image of $g$ is not conjugate to any element of the image of $A$, as needed.
\end{pf}

This confirms hypothesis 2 of Theorem \ref{t:Combination theorem}.

\begin{lem}\label{edge intersection}
Let $M$ be a large Seifert-fibred manifold and let $C_1, C_2$ be distinct peripheral subgroups of $G=\pi_1(M)$.  Then $\widehat C_1\cap \widehat C_2=\widehat Z$, the profinite completion of the canonical normal cyclic subgroup of $G$.
\end{lem}
\begin{pf}
Let $p$ be the epimorphism from $\pi_1(M)$ to $\pi_1(O)$.  Since $\pi_1(O)$  is virtually free,  the intersection of $\widehat{p(C_1)}$ and $\widehat{p(C_2)}$ is trivial since the intersection of $p(C_1)\cap p(C_2)$ is (see Lemma 3.6 in \cite{ribes_conjugacy_1996}), so the intersection of $\widehat C_1$ and $\widehat C_2$ is $\widehat Z$.
\end{pf}

This confirms hypothesis 3 of Theorem \ref{t:Combination theorem}.  Lastly, we need to prove that the JSJ decomposition is profinitely 2-acylindrical.

\begin{lem} \label{Profinitely acylindrical}
Let $M$ be a graph manifold in which every Seifert-fibred piece is large, let $G=\pi_1(M)$ and let $(\G,\Gamma)$ be the JSJ decomposition of $G$.  Then $(\G,\Gamma)$ is profinitely 2-acylindrical.
\end{lem}
\begin{pf}
By Lemma \ref{l: 2-acylindrical lemma}, $G$ is 2-acylindrical.   Choose three consecutive edges $e_1,e_2,e_3$ in $S(\widehat G)$. We have to prove that the intersection $\widehat G_{e_1}\cap \widehat G_{e_2}\cap \widehat G_{e_3}$ of their edge stabilizers  is trivial.  By translating them if necessary we may assume that $e_2$ is in $S=S(G)$ and so its common vertices $v,w$ with $e_1$ and $e_3$ respectively are in $S$ as well.   Let $Z_v$ and $Z_w$ be the canonical normal cyclic subgroups of $G_v$ and $G_w$ respectively.  By Lemma \ref{edge intersection}, we have $\widehat G_{e_1}\cap \widehat G_{e_2}=\widehat Z_v$ and $\widehat G_{e_2}\cap \widehat G_{e_3}=\widehat Z_w$.  By the minimality of the JSJ decomposition of $M$, $Z_v$ and $Z_w$ have trivial intersection and hence $Z_v\times Z_w$ is a finite-index subgroup of $G_{e_2}$.   But $\widehat{Z_v\times Z_w}=\widehat Z_v\times \widehat Z_w$, so we have
\[
\widehat G_{e_1}\cap \widehat G_{e_2}\cap \widehat G_{e_3}=\widehat Z_v\cap\widehat Z_w=1
\]
as required.
\end{pf}

We can now prove Theorem \ref{t: Graph manifolds are c.s.}.

\begin{pfof}{Theorem \ref{t: Graph manifolds are c.s.}}
Let $M$ be an orientable graph manifold.  After doubling along the boundary, we can assume that $M$ is closed.  If $M$ is Seifert-fibred then its fundamental group is conjugacy separable by Theorem \ref{t: Martino}.  Otherwise, the JSJ decomposition of $M$ is non-trivial and, by Theorem \ref{t:Efficiency}, the profinite topology on $G=\pi_1(M)$ equipped with the JSJ decomposition is efficient.

If $M$ is finitely covered by a torus bundle over a circle then either $M$ is Seifert-fibred or $M$ admits a geometric structure modelled on the Lie group $\mathrm{Sol}$.  In this second case, it is a standard fact that $M$ has a normal cover $M'$ of degree at most 8 such that $M'$ is a torus bundle over the circle (see, for instance, Theorem 4.17 of \cite{scott_geometries_1983}).   Because all groups of order at most 8 are polycyclic and $\pi_1(M')$ is polycyclic, it follows that $G$ is polycyclic. Therefore $G$ is conjugacy separable by the main theorem of \cite{remeslennikov_conjugacy_1969}.

Now we assume that $M$ is not finitely covered by a torus bundle, so some pieces of the JSJ decomposition are large.  We first deal with the special case in which every piece of the JSJ decomposition $\G$ of $G$ is large.   Then $\G$ is profinitely 2-acylindrical by Lemma \ref{Profinitely acylindrical}.  Lemmas \ref{conjugacy distinguished} and \ref{edge intersection}, together with Niblo's theorem that the fundamental groups of Seifert-fibred manifolds are double-coset separable \cite{niblo_separability_1992}, imply that the hypotheses of Theorem \ref{t:Combination theorem} hold.

Now suppose that some of the pieces of $M$ are not large.  As in the proof of Lemma \ref{l: 4-acylindrical lemma}, $G=\pi_1(M)$ has an index-two subgroup $G'$ corresponding to a two-sheeted cover $M'$ such that every piece of the JSJ decomposition of $M'$ is large.  Therefore, by the preceding paragraph, $G'$ is conjugacy separable.  Fix a coset representative $t\in G\smallsetminus G'$.

Let $a,b\in G$ and suppose that $b$ is in the profinite closure of the conjugacy class $a^G$.  We need to prove that $b$ is conjugate to $a$.  Suppose first that $a\in G'$.  The conjugacy class $a^{G'}$ is closed in $G'$ and hence in $G$; likewise, $(a^t)^{G'}$ is closed in $G'$, so $a^G=a^{G'}\cup (t^{-1}at)^{G'}$ is closed in $G$, so $b\in a^G$ as required.

It remains to consider the case when $a\in G\smallsetminus G'$ and hence $b\in G\smallsetminus G'$ also.  As $b\in \overline{a^G}$ it follows that $b^2$ is in the closure of the conjugacy class of $a^2$ and hence, by the previous paragraph, $b^2$ is conjugate to $a^2$.  Replacing $b$ with a conjugate, therefore, we may assume that $a^2=b^2$.

Consider the action of $a^2$ on the Bass--Serre tree $T$ of $\G$.  If $a^2$ is hyperbolic then $a$ and $b$ are both hyperbolic, and we have that
\[
T_a=T_{a^2}=T_{b^2}=T_b
\]
and moreover $a$ and $b$ have the same translation length.  Therefore $ab^{-1}$ fixes $T_a$ pointwise and so, by Lemma \ref{l: 4-acylindrical lemma}, $a=b$ as required.

If $a^2$ is elliptic then $a$ is also elliptic and stabilizes some vertex $u'\in T$.  Because $a\notin G'$, the stabilizer $G_u$ of $u'$ is isomorphic to $\Z\rtimes\Z$ and $a$ is not contained in the stabilizer of the two incident edges, which both equal the abelian subgroup of index two and are contained in $G'$.  (To see this, note that in the construction of $M'$ in the proof of Lemma \ref{l: 4-acylindrical lemma}, every torus of the JSJ decomposition of $M$ lifts to $M'$.)  The same argument shows that $b$ stabilizes a vertex $v'$ in a similar fashion.

Suppose $u'\neq v'$.  The square $a^2=b^2$ stabilizes both vertices.  As $M$ is not finitely covered by a torus bundle, $u'$ and $v'$ are not adjacent, so $a^2$ also stabilizes a vertex $w'\in T$, adjacent to $u'$, that corresponds to a large piece of the JSJ decomposition.  Furthermore, $a^2$ stabilizes two edges incident at $w'$, so $a^2$ is contained in the normal cyclic subgroup of $G_w$.  But $a^2$ is also contained in a normal cyclic subgroup of $G_u$, which contradicts the minimality of the JSJ decomposition of $M$.

Therefore,  $a,b\in G_u$ and $u'$ is the unique vertex stabilized by $a$ and $b$.  For any $g\in G$, $gu'$ is the unique vertex stabilized by $gag^{-1}$ and so if $gag^{-1}\in G_u$ then $g\in G_u$.  That is, $a^G\cap G_u=a^{G_u}$ and so
\[
b\in \overline{a^G}\cap G_u=\overline{a^G}\cap \overline{G}_u\cap G_u=\overline{a^G\cap G_u}\cap G_u=\overline{a^{G_u}}\cap G_u.
\]
Theorem \ref{t:Efficiency} implies that the profinite topology of $G$ induces the full profinite topology on $G_u$, so $b$ is in the closure of $a^{G_u}$ in the profinite topology on $G_u$.  As $G_u$ is conjugacy separable, $b$ is conjugate to $a$ in $G_u$ as required.
\end{pfof}

\bibliographystyle{plain}

\bigskip\bigskip\centerline{\textbf{Authors' addresses}}
\smallskip\begin{center}\begin{tabular}{ll}%
Henry Wilton & Pavel Zalesskii\\
Department of Mathematics & Department of Mathematics\\
1 University Station C1200 & University of Brasilia\\
Austin, TX 78712-0257 & 70910-900 Brasilia-DF\\
USA & Brazil\\
{\texttt{henry.wilton@math.utexas.edu}}&{\texttt{pz@mat.unb.br}}
\end{tabular}\end{center}

\end{document}